\documentclass{amsproc}      
\usepackage{amssymb}  
\usepackage{amsmath}  
\usepackage{amsthm}  
\newcommand{\be}[1]{\begin{eqnarray#1}}
\newcommand{\ee}[1]{\end{eqnarray#1}} 
\newtheorem{thm}{Theorem}[section]
\newtheorem{propn}[thm]{Proposition}      
\newtheorem{lemma}[thm]{Lemma}

\theoremstyle{remark}
\newtheorem{remark}[thm]{Remark}
\newtheorem{definition}[thm]{Definition}

\newcommand{\tp}{\otimes}

\newcommand{\braid}[2]{{#1}$\lower4pt\hbox{$\tp\atop\raise4pt
            \hbox{$\scriptscriptstyle\Ru $}$}${#2}}
\newcommand{\twist}[2]{{#1}${\,\scriptscriptstyle \Ru}\atop\raise9pt\hbox{$\scriptstyle\tp$} ${#2}}
\newcommand{\twistF}[2]{{#1}${\,\scriptscriptstyle \F}\atop\raise9pt\hbox{$\scriptstyle\tp$} ${#2}}

\newcommand{\Ru}{\mathcal{R}}

\newcommand{\U}{\mathcal{U}}                
\newcommand{\F}{\mathcal{F}}
\newcommand{\C}{\mathbb{C}}                  
\newcommand{\Z}{\mathbb{Z}}                  
                  
\newcommand{\N}{\mathbb{N}}

\newcommand{\la}{\lambda}

\newcommand{\n}{\nonumber          }           
\newcommand{\al}{\alpha}                         
\newcommand{\bt}{\beta}

\newcommand{\card}{\mathrm{card}}              
\newcommand{\End}{\mathrm{End}}

\newcommand{\T}{\mathrm{T}}
\newcommand{\si}{\sigma}

\begin{document}
\title{Characters of $\U_q\bigl(gl(n)\bigr)$-reflection equation algebra}
\author{A. Mudrov}
\address{Department of Mathematics, Bar-Ilan University, 52900 Ramat-Gan, Israel}
\email{mudrova@macs.biu.ac.il}
\maketitle
\begin{abstract}
We list characters (one-dimensional representations) of the reflection
equation algebra associated with the fundamental vector representation 
of the Drinfeld-Jimbo quantum group $\U_q\bigl(gl(n)\bigr)$.
\end{abstract}
\section{Introduction}
Let $V$ be a complex vector space and $S\in \End^{\tp2}(V)$
a Yang-Baxter operator, satisfying the braid identity
$$
S_{12}S_{23}S_{12}=S_{23}S_{12}S_{23}.
$$
Here, $S_{12}=S\tp 1$ and $S_{23}=1\tp S$.
A matrix $A\in \End(V)$ is called a solution to the numerical 
reflection equation (RE) or a numerical RE matrix if the 
equality 
\be{}
\label{re}
SA_2 S A_2 = A_2 S A_2S,
\ee{}
where $A_2=1\tp A$,
holds in $\End^{\tp2}(V)$. Numerical RE matrices define 
one-dimensional representations (characters) of the RE algebra,
the quotient of the tensor algebra $\T\bigl(\End^*(V)\bigr)$
by quadratic relations (\ref{re}), \cite{KSkl,KS}.
RE and related algebraic structures appeared in the theory of
integrable models \cite{Cher,Skl,AFS}; 
they find applications in 3-dimensional topology \cite{K}  and problems 
of  covariant quantization on
G-spaces \cite{Do1,Do2}. In the latter case, they are related to 
representations 
of quantum groups on $V$; then the operator $S$ is equal to $PR$,
the product of the permutation $P$ in $V\tp V$ and the image
$R$ of the universal
R-matrix, \cite{Dr}. Numerical solutions to equation 
(\ref{re}) are the key ingredient in 
the character method of covariant quantization on homogeneous spaces
developed in \cite{DoM1,DoM2}. 

In the present paper we study equation (\ref{re}) associated
with the fundamental vector representation of the quantum
group $\U_q\bigl(gl(n)\bigr)$. 
An example of solution to this equation 
is the matrix 
$D_n= \la \sum^n_{i=1}  e^i_{n+1-i}$,
$\la\in \C$, 
as well as its similarity transformation by a diagonal
invertible matrix, according to \cite{KSS}.
Therein, all the non-degenerate RE matrices 
were presented in dimensions $n=2,3,4$:
$$
A^{1,1}=
\left( 
\begin{array}{cc}
 \la+\mu  & y_1  \\
 y_2 & 
\end{array} 
\right) ,
\quad A^{2,1}=
\left(  
\begin{array}{ccc}
 \la+\mu & & y_1 \\ 
& \la &  \\  
y_3  &  &  
\end{array}  
\right),
$$  
$$
A^{2,2}= 
\left( 
\begin{array}{cccc}
 \la+\mu & & & y_1 \\
 & \la+\mu &y_2& \\ 
 &y_3 & & \\ 
y_4 &   &   &    
\end{array}   
\right)    ,\quad    
A^{3,1}= 
\left(
\begin{array}{cccc}
 \la+\mu& & & y_1 \\
 & \la & & \\
 & & \la & \\
 y_4 & & &
\end{array}
 \right),
$$ 
$$
\mbox{where}\quad y_iy_{n+1-i}=-\la\mu\not=0.
$$
Here, the parameterization is adapted to our exposition.
Besides, in \cite{DoM1}, there were found diagonal numerical 
RE matrices $P_k=\la \sum_{i=1}^k e_i^i$, $k\leq n$.
The above mentioned matrices exhaust all (as far as we can judge) 
known solutions
to (\ref{re}) associated to the standard quantum group
$\U_q\bigl(gl(n)\bigr)$. In the present paper, we list all
possible solutions.
\section{Result}
The image of the universal R-matrix of
the Drinfeld-Jimbo quantum group $\U_q\bigl(gl(n)\bigr)$
in the fundamental vector representation on $\C^n$
is, \cite{FRT},
\be{}
\label{rsln}
R &=& q\sum_{i=1,\ldots, n} e^i_i \tp e^i_i +
\sum_{i,j=1,\ldots, n \atop i\not =j} e^i_i \tp e^j_j +
\omega \sum_{i,k=1,\ldots, n \atop i < k} e^k_i \tp e^i_k.
\ee{}
Here, $\{e^i_j\}$ is the standard basis in $\End(\C^n)$ with the multiplication 
$e^i_j e^l_k =  \delta^l_j  e^i_k $ expressed through
the Kronecker symbols $\delta^l_j $; 
$\omega$ stands for $=q-q^{-1}$.
It is convenient to represent the corresponding  braid  matrix
in the form 
\be{}
\label{S}
S &=& \sum_{i,k=1,\ldots, n } s_{ik} e^i_i \tp e^k_k +
\sum_{i,j=1,\ldots, n \atop i\not =j} e^i_j \tp e^j_i,
\>\mbox{where}\;
s_{ik}=
\left\{
\begin{array}{ll}
\omega,& i <k,\\
q,& i =k,\\
0,& i >k,
\end{array}
\right.
\ee{}
We shall deal with subsets in $\Z$ and denote by $[a,b]$ the intervals
$\{k\in \Z| a\leq k \leq b\}$. Instead of the square
brackets, we use parentheses for the
intervals defined by strict inequalities.
To formulate the classification theorem,  we need the following data.
\noindent
\begin{definition}
\label{ap}
{\em An admissible  pair} $(Y,\si)$ consists of an
ordered subset 
$Y\subset I=[1,n]$  and a decreasing injective map 
$\si\colon Y \to I$ without stable points.
\end{definition}
\noindent
Clearly, the map $\si$ is determined by its image $\si(Y)$.
 In the set $Y$, we distinguish two
non-intersecting subsets $Y_+=\{i\in Y| i>\si(i)\}$ 
and $Y_-=\{i\in Y| i<\si(i)\}$; obviously $Y=Y_-\cup Y_+$. Denote
$b_-=\max\{Y_-\cup \si (Y_+)\}$ and  $b_+=\min\{Y_+\cup \si (Y_-)\}$.
Because the injection $\si$ is decreasing, $b_-<b_+$.
We adopt the convention $b_-=0$ and $b_+=n+1$ if
$Y=\emptyset$.
\begin{thm}[Classification]
\label{clth}
General  solution to equation (\ref{re})
with the braid matrix (\ref{S})
has the form
\be{}
\label{gf}
A&=&\sum_{i\in I} x_i e^i_i + \sum_{j\in Y}y_j  e^{\si(j)}_j,
\hspace{28pt}
\ee{}
where 
\begin{itemize}
\item $(Y,\si)$ is an  admissible pair,
\be{}
Y&=&[1, b_-]\cup [b_+, b_+ + b_- -1] ,\hspace{20pt}
\\
\si(i)&=&b_+ + b_- -i
\ee{}
for $b_-,b_+\in \N$ subject to the conditions
$b_-<b_+$, $b_- + b_+\leq n+1$ and 
\be{}
x_i & = &
\left\{
\begin{array}{ll}
\la+\mu,& i\in [1,b_-],\\
\la,& i\in (b_-,b_+),\\
0,& i\in [b_+,n],
\end{array}
\right.\hspace{35pt}\\
y_i y_{\si(i)} & = &-\la\mu \not =0
\ee{}
for  $\la,\mu \in \C$;
\item  $(Y,\si)$ is an admissible pair such that
 $Y \cap \si(Y)=\emptyset$ and
\be{}
x_i & = &
\left\{
\begin{array}{cc}
\la,& i\in [1,b],\\
0,& i\in (b,n],
\end{array}
\right.\hspace{42pt}\\
y_i &  \not= &0
\ee{}
for $b\in [b_-,b_+)$ and $\la\in \C$.
\end{itemize}
\end{thm}
\begin{remark}
As follows from the theorem, there are two classes of 
numerical RE matrices, those corresponding to 
$\si(Y)=Y$ and $\si(Y)\cap Y=\emptyset$. We call them
solutions of Type 1 and 2, respectively.
\end{remark}
\section{Proof of the classification theorem}
The proof of Theorem (\ref{clth}) formulated in the previous section 
is combinatorial. It is a result of the direct analysis of 
equation (\ref{re}) organized into a sequence of lemmas.
Putting $A=\sum_{\al,\bt=1}^n A^\al_\bt e^\bt_\al$, we compute
$$
A_2 S A_2 = 
\sum_{i,\bt,\al,\nu} s_{i\nu}  A^{\nu}_\bt A^\al_{\nu} \; e^i_i\tp e^\bt_\al
 +        
\sum_{i,j,\bt,\al\atop i\not=j} A^{j}_\al A^\bt_{i} \; e^i_j\tp e^\al_\bt .
$$
Substituting this into the left- and right-hand sides of (\ref{re}),
we rewrite it in a more explicit form:
\be{}
SA_2SA_2&=& 
\sum_{i,\bt,\al,\nu} s_{i\bt}s_{i\nu} A^{\nu}_\bt A^\al_{\nu} \; e^i_i\tp e^\bt_\al 
+
\sum_{i,j,\bt,\al \atop j\not = i} s_{i\al}A^{j}_\al A^\bt_i   \; e^i_j\tp e^\al_\bt 
\label{lhs}\\
&+& 
\sum_{i,\bt,\al,\nu \atop i\not = \bt} s_{i\nu} A^{\nu}_\bt A^\al_{\nu} \; e^\bt_i\tp e^i_\al 
+
\sum_{i,j,\bt,\al \atop {j\not = i \not = \al }} A^{j}_\al A^\bt_i   \; e^\al_j\tp e^i_\bt 
,\n\\
A_2SA_2S&=& 
\sum_{i,\bt,\al,\nu} s_{i\al}s_{i\nu} A^{\nu}_\bt A^\al_{\nu} \; e^i_i\tp e^\bt_\al 
+
\sum_{i,j,\bt,\al \atop j\not = i} s_{j\bt}A^{j}_\al A^\bt_i   \; e^i_j\tp e^\al_\bt 
\label{rhs}\\
&+& 
\sum_{i,\bt,\al,\nu \atop i\not = \al} s_{i\nu} A^{\nu}_\bt A^\al_{\nu} \; e^i_\al\tp e^\bt_i 
+
\sum_{i,j,\bt,\al \atop {i\not = j \not = \bt }} A^{j}_\al A^\bt_i   \; e^i_\bt\tp e^\al_j
\n.
\ee{}
Comparison of (\ref{lhs}) with (\ref{rhs}) gives rise to the system
of quadratic equations on the matrix elements $A^i_j$.
\begin{lemma}
Equation (\ref{re}) is equivalent to the following system of equations:
\be{}
\left\{
\begin{array}{rrl}
A^m_i A^n_i &=&0, 
 \\
A^i_m A^i_n &=&0, %& m\not = n  \not = i \not = m
\end{array}
\right.
\hspace{44pt}
\quad  m\not = n  \not = i \not = m ,
\hspace{20pt}
\label{eq1}
\\
A^n_i A^j_m \hspace{8pt}=\hspace{8pt}0, 
\hspace{36pt}
\quad
\left\{ 
\begin{array}{l}
j\not = m \not = n\not =i,\\
(m-i)(n-j)<0,
\end{array}
\right.
\label{eq2}
\\
\left\{
\begin{array}{rrll}
(q - s_{im}) A^i_i A^m_i &=& \sum_{\nu}s_{i\nu} A^{\nu}_i A^m_{\nu},& i \not = m
, \\
(q - s_{im}) A^i_i A^i_m &=& \sum_{\nu}s_{i\nu} A^{\nu}_m A^i_{\nu},& i \not = m
, \\
0&=&\sum_{\nu}s_{i\nu} A^{\nu}_m A^n_{\nu} ,&
(m-i)(n-i)<0,
\end{array}
\right.
 \label{eq3}\\
A^n_i A^i_m - \sum_{\nu}s_{i\nu} A^{\nu}_m A^n_{\nu}=
(s_{ni} - s_{im}) A^i_i A^n_m 
, \quad  m\not = i  \not = n \not = m
, \hspace{8pt}
\label{eq4}
\\
\omega A^m_m A^i_i
=
\sum_{\nu}s_{i\nu} A^{\nu}_m A^m_{\nu}  -  \sum_{\nu}s_{m\nu}  A^{\nu}_i
A^i_{\nu}
, \quad  i< m
%\hspace{pt}
\label{eq5}
\ee{}
\end{lemma}
\begin{proof}
This statement is verified by the direct analysis.
\end{proof}
The next lemma accounts for equations  (\ref{eq1})
and  (\ref{eq2}). 
\begin{lemma}
\label{l-eq1}
If $A$ a solution to equation (\ref{re}), then it
can be represented in the form (\ref{gf}),
%$$A=\sum_{i\in I} x_j e^j_j + \sum_{i\in Y}y_i e^i_{\si(i)},$$
where $(Y,\si)$ is an admissible pair, $x_i=A^i_i$  for $i\in I$,
and  $y_i=A^{\si(i)}_i\not = 0$ for $i\in Y$.
\end{lemma}
\begin{proof}
By virtue of equations (\ref{eq1}), the matrix $A$ has 
at most one non-zero off-diagonal entry in 
every row and every column. 
Indices of such rows  form a subset $Y$ in $I$,
and the  off-diagonal entries may be written as 
$y_i=A_i^{\si(i)}\not = 0$,  $i\in Y$, for some 
bijection $\si\colon Y\to I$ with no stable points. 
This is equivalent to equation (\ref{eq1}). The
map $\si$ is decreasing; that is encoded in equation 
(\ref{eq2}). 
\end{proof}
\begin{lemma}
\label{l-eq31}
The equations of system (\ref{eq3}) 
are equivalent to
\be{}
\label{eq3'}
\sum_{\al\geq\max(i,m)}A^\al_i A^m_\al=0,\quad i\not = m,
\ee{}
by virtue of (\ref{eq2}). They imply $x_i=0$ for  $i = b_+$.
\end{lemma}
\begin{proof}
Setting $i=j$ in (\ref{eq2}) and substituting it into 
the equations of system (\ref{eq3}), we reduce them
to (\ref{eq3'}). 
Setting  $m =\si(i) < i \in Y_+$ in (\ref{eq3'})
leads to $ x_i y_i =0$ and therefore $x_i=0$.
Similarly, the assumption  $m=\si(i)>i \in Y_-$ reduces 
(\ref{eq3'}) to $ y_i x_{\si(i)}=0$ and hence $x_{\si(i)}=0$. 
So $x_i=0$ for all $i\in Y_+\cup \si(Y_-)$
and $i=b_+$ in particular
\end{proof}
\begin{lemma}
\label{l-eq4}
Lemma \ref{l-eq31} taken into account, 
equation (\ref{eq4}) is equivalent to the following two assertions.
\begin{enumerate}
\item For any $m\in Y$, either $\si^2(m)=m$ or  $\si(m)\not\in Y$.
\item $x_i = x_{b_-} $ and $x_j = 0$ whenever $i\leq b_-$ and $j\geq b_+ $.
\end{enumerate}
\end{lemma}
\begin{proof}
First note that if $m\not \in Y$, equation (\ref{eq4})
holds identically.
So we can assume $m \in Y$ and rewrite (\ref{eq4}) as
\be{}
A^n_i A^i_m - s_{im} A^{m}_m A^n_{m}-
s_{i\si(m)} A^{\si(m)}_m A^n_{\si(m)}=   
(s_{ni} - s_{im}) A^i_i A^n_m,
\label{eq4'}
\ee{}
$$
m\not=i\not=n \not=m.
$$
Supposing $n\not =\si(m)$ we find 
$A^n_i A^i_m - s_{i\si(m)} A^{\si(m)}_m A^n_{\si(m)}= 0.$ 
If $i=\si(m)$, then, having in mind
$s_{ii}=q\not=1$ and $A^{\si(m)}_m=y_m\not =0$, we obtain
$A^n_{\si(m)}=0$. Since  $n\not =\si(m)=i$ and $n\not =m$,
this means 
either ${\si(m)}\not \in Y$ or ${\si^2(m)}=m$. Assuming
$i\not=\si(m)$, we come to the equation 
$s_{i\si(m)} A^{\si(m)}_m A^n_{\si(m)}=0$, which is fulfilled
as well, due to $A^n_{\si(m)}=0$.

It remains to study the case  $n =\si(m)$.
Under this hypothesis, the term $A^n_i A^i_m$ vanishes. 
Indeed, since $i\not= m$, 
one has $A^i_m \not =0 \Rightarrow i= \si(m)$. But this
contradicts the condition $i\not= n =\si(m)$.
In terms of the variables $x_i$ and $y_i$, equation (\ref{eq4'})
reads
$$
- s_{im} x_m y_m-s_{i\si(m)} y_m x_{\si(m)}=   
(s_{\si(m)i} - s_{im}) x_i y_m,
$$
$$
m\not=i\not =\si(m).
$$
Depending on allocation of the indices $i$, $m$, and $\si(m)$, 
this equation splits into the following four implications.
\be{}
i<m\quad\mbox{and} \quad i<\si(m) &\Longrightarrow&  x_i=x_m+x_{\si(m)}
,\label{eq41}\\
i<m\quad\mbox{and} \quad i>\si(m) &\Longrightarrow&  x_m=0
\label{eq42},\\
i>m\quad\mbox{and} \quad i<\si(m) &\Longrightarrow&  x_{\si(m)}=0
\label{eq43},\\
i>m\quad\mbox{and} \quad i>\si(m) &\Longrightarrow&  x_i=0
\label{eq44}.
\ee{}
Recall that the index $m$ is assumed to be from $Y$.
Equation (\ref{eq44}) is equivalent to $x_i = 0$ for $i>b_+$.
By Lemma \ref{l-eq31}, this is also true for $i\geq b_+$.
Equations (\ref{eq42}) and  (\ref{eq43}) are thus satisfied
as well.
Equation (\ref{eq41}) states 
$x_i = x_{b_-}+x_{\si(b_-)}$
if $i<b_- \in Y_-$
or 
$x_i = x_{b_-}+x_{\si^{-1}(b_-)}$
if $i<b_-\in \si(Y_+)$.
Applying Lemma \ref{l-eq31}, we find 
$x_i = x_{b_-}$, in either cases.
\end{proof}
\begin{remark}
We used Lemma \ref{l-eq31} in the proof of Lemma \ref{l-eq4}
in order to find the values  $x_{b_\pm}$ on the boundaries
of the intervals $[1,b_-]$ and  $[b_+,n]$. 
If $(b_-,b_+)\not=\emptyset$,
equations (\ref{eq42}) and  (\ref{eq43}) are enough for
that purpose, because we can put $i\in (b_-,b_+)$ there and
come to the same result.
\end{remark}
\begin{lemma}
\label{l-eq3}
Equation (\ref{eq3'}) is fulfilled 
by virtue of Lemma \ref{l-eq4}.
\end{lemma}
\begin{proof}
Consider the case $m<i$. Equation  (\ref{eq3'}) holds
if $i\not \in Y$, because the sum turns into $A^i_iA^m_i$. 
So we may assume $i \in Y$ and distinguish two
cases: $i\in Y_-$ and $i\in Y_+$. 
Assumption $i<\si(i)$ leads to
$A^i_i A^m_i + A^{\si(i)}_i A^m_{\si(i)}=0$.
The equality $m={\si(i)}$ is impossible, since otherwise
$m={\si(i)}<i<{\si(i)}$. Therefore the first term
vanishes and we come to $A^m_{\si(i)}=0$, $m\not= i$.
This condition is satisfied, by Lemma  \ref{l-eq4},
Statement 1.
The case  $i\in Y_+$ results in 
$x_i A^m_i=0$. Setting $m =\si(i)$, we come to $ x_i =0 $.
Once $i\in Y_+\Rightarrow i\geq b_+$, we 
encounter a particular case of Lemma \ref{l-eq4},
Statement 2.

We should study the situation $i<m$. Equation
(\ref{eq3'}) evidently holds if $\si(i)<m$. 
We may think that  $m\leq\si(i)$; then
$$
\begin{array}{llrll}
i<m=\si(i)& \Rightarrow & A^m_m = 0& \Rightarrow & x_{\si(i)}=0,\\
i<m<\si(i)& \Rightarrow &  A^m_{\si(i)} = 0& \Rightarrow &
\si(i)\not \in Y.
\end{array}
$$
These requirements are fulfilled, by Lemma  \ref{l-eq4}.
\end{proof}

It remains to satisfy equation (\ref{eq5}), to complete the 
proof of Theorem \ref{clth}.
Let $Y_0$ be the subset in $Y$, such that $\si$ restricted to $Y_0$ is
involutive. By Lemma \ref{l-eq4}, either $\si(i) \not \in Y$
or $i$ and $\si(i)$ belong to $Y_0$ simultaneously.
Equation (\ref{eq5}) falls into the four equations
\be{}
       x_m(x_i-x_m) & = & 0 , \hspace{134pt}  i\not \in Y_0,\> m\not \in Y_0,
\label{eq5'}\\
\omega x_m(x_i-x_m) & = & s_{i\si(m)} y_m y_{\si(m)},\quad 
                  \hspace{66pt} i \not\in Y_0,\> m \in Y_0,
\label{eq5''}\\
\omega x_m(x_i-x_m) & =  &-s_{m\si(i)} y_i y_{\si(i)} ,\quad
                  \hspace{67pt} i \in Y_0,\> m \not \in Y_0,
\label{eq5'''}\\
\omega x_m(x_i-x_m) & =  & s_{i\si(m)} y_m y_{\si(m)}
- s_{m\si(i)} y_i y_{\si(i)} ,\quad
        i \in Y_0,\> m \in Y_0.
\label{eq5''''}
\ee{}
Everywhere $i<m$, see (\ref{eq5}).
\begin{lemma}
\label{l-eq51}
Suppose $Y_0\not=\emptyset$. 
Then $Y=Y_0$ and, moreover,
\be{}
Y_-=\{1,\ldots,b_-\} ,\quad Y_+=\{b_+,\ldots,b_+ +b_- -1\},
\label{notempty}
\ee{}
\be{}
\si(i) = b_+ + b_- -i.
\ee{}
\end{lemma}
\begin{proof}
Let $k \in Y_0 \cap Y_-$ and  $\si(k)\in Y_0 \cap Y_+$ .
Suppose either $l \in Y_-\backslash Y_0$ 
or $l<\min(Y_-)$. 
The assumption
$l<k$ contradicts  equation (\ref{eq5''})
if one sets $i=l$, $m=k$. The 
inequality $k<l$ does not agree with (\ref{eq5'''})
if one sets $i=k$, $m=l$.
In both cases one uses Lemma \ref{l-eq4}, Statement 2, and gets
$0=y_k y_{\si(k)}$; that is impossible since
$y_k\not=0$ for all $k\in Y$, by definition of $Y$.
Therefore $Y_-\subset Y_0$ and 
$\min(Y_-)=1$. 
Assume now  $l \in Y_+\backslash Y_0$. Then,
$l$ cannot exceed $\si(k)$,
because otherwise 
$\si(l) <\si^2(k) = k \Rightarrow\si(l)\in Y_0\Rightarrow l\in Y_0$,
the absurdity.
The only possibility is $l<\si(k)$. But this again contradicts 
equation  (\ref{eq5'''}) if one sets $i=k$ and $m=l$.
\end{proof}
\begin{lemma}
\label{l-eq52}
Suppose $Y_0\not=\emptyset$. 
Then, there is $a\in \C$ such that $y_i y_{\si(i)} = a$ for all $i\in Y$. 
If $(b_-,b_+)\not=\emptyset$, then   $x_m=\la$ for all $m \in (b_-,b_+)$,
where $\la$ is a solution 
to the quadratic equation $\la(\la-x_{b_-})=a$.
\end{lemma}
\begin{proof}
By Lemmas \ref{l-eq4} and \ref{l-eq51}, 
equation (\ref{eq5''''}) is fulfilled if and only if
the product $y_i y_{\si(i)}$ does not depend on $i\in Y$
and is equal to some $a \in \C$.
If  $(b_-,b_+)\not=\emptyset$, equation  (\ref{eq5'''})
suggests $x_m(x_m-x_{b_-})=a$ as soon as $m\in (b_-,b_+)$.
It is easy to see that equation  (\ref{eq5'}) is equivalent
to $x_m=x_i=\la$
for $i,m\in (b_-,b_+)$ and some $\la\in \C$. 
\end{proof}
It is convenient to introduce the parameterization  $x_{b_-}=\la+\mu$, 
$a=-\mu\la$.
Then $x_m = \la + \mu$ for $m\leq b_-$, 
$x_m= \la$ for  $b_-<m<b_+$, and $y_i y_{\si(i)}=-\la\mu$.
This parameterization makes sense even if  $(b_-,b_+)=\emptyset$.
The scalars $\la$ and $\mu$ will have  the meaning of 
eigenvalues of the matrix $A$.
\begin{lemma}
Let $A$ be a matrix satisfying equations (\ref{eq1})-(\ref{eq4}) 
and $(Y,\si)$ the corresponding admissible pair.
Then, equation (\ref{eq5}) gives rise to the 
alternative 
\begin{itemize}
\item
$Y_0=Y$. $A$ is a solution of Type 1 from  Theorem \ref{clth}.
\item
$Y_0=\emptyset$. 
$A$ is a solution of Type 2 from Theorem \ref{clth}.
\end{itemize}
\end{lemma}
\begin{proof}
If $Y_0=\emptyset$, then $Y\cap \si(Y)=\emptyset$ by 
Lemma \ref{l-eq4}, Statement 1. Equation 
(\ref{eq5}) is reduced to (\ref{eq5'}). It is satisfied if and only if 
$x_m=\la$, $i\leq b$, and $x_m=0$, $b< i$, for 
some $\la\in \C$ and  $b\in \Z$. Comparing this
with Lemma \ref{l-eq4}, Statement 2,
we come to a solution of Type 2 with  $b\in [b_-,b_+)$.

Suppose $Y_0\not= \emptyset$. Then, 
equation (\ref{eq5}) implies Lemmas \ref{l-eq51}
and  \ref{l-eq52}. Conversely, these lemmas
ensure (\ref{eq5'})--(\ref{eq5''''}), which are 
equivalent to  (\ref{eq5}). Altogether, this
gives a solution of Type 1.
\end{proof}
Thus we accomplish the proof of Theorem \ref{clth} and proceed
to the analysis of the solutions obtained.
\section{Structure of solutions}
The first question about the structure of the numerical RE matrices 
listed in Theorem \ref{clth} is
what pairs $(Y,\si)$ may participate in the classification.
The solutions of Type 1 involve two integers $b_-$, $b_+$ subject 
to the conditions of the theorem. They determine the 
admissible pair $(Y,\si)$.
The solutions of Type 2 are labeled with subsets $Y\subset I$ and 
injective maps $\si\colon Y\to I\backslash Y$. Clearly, 
$K=\card(Y)\leq\frac{n}{2}$ and, with $Y$ given, there are $C_{n-K}^K$ 
possibilities for $\si$, i.e. the number
of subsets in $I\backslash Y$ with $\card=K$.
Further we describe the properties of solutions to 
equation (\ref{re}).
We consider the standard basis  $\{e^i\}$  in $\C^n$;
the action of $\End(\C^n)$ on $\C^n$ being  given
by $e^i_j e^k = \delta^k_j e^i$. 

\begin{propn}
Let $(Y,\si)$ be an admissible pair corresponding
to a numerical RE matrix $A$.
The subspaces 
$V^2_i=\C e^i\oplus \C e^{\si(i)}$, $i\in Y$, and 
 $V^1_i=\C e^i$, $i \not\in Y\cup \si(Y)$, are $A$-invariant.
\end{propn}
\begin{proof}
This follows from the general form (\ref{gf}) of the RE matrices 
and from Lemma \ref{l-eq4}, Statement 1.
\end{proof}
The spectral properties of numerical RE matrices are described by
the following proposition.
\begin{propn}
An RE matrix  of Type 1 from Theorem \ref{clth} has 
eigenvalues $\mu$, $\la$, and $0$ of 
multiplicities $b_-$, $b_+ - b_- -1$, and 
$n- (b_- + b_+) + 1$, respectively.
It is semisimple
if and only if $\la\not = \mu$.

An RE matrix  of Type 2 from Theorem \ref{clth} has 
eigenvalues $\la$ and $0$ of 
multiplicities $b$ and $b-K$
respectively.
It is semisimple
if and only if $\la\not = 0$.
\end{propn}
\begin{proof}
Restricted to $V^2_i$, $i\in [1,b_-]$, an RE matrix
$A$  of  of Type 1 is equal to 
$
\begin{array}{||cc||}
\la+\mu& y_i\\
y_{\si(i)}&0
\end{array}\:
$. Because $y_iy_{\si(i)}=-\la\mu$, this $2\times2$ matrix
has eigenvalues $\mu$ and $\la$. 
On the subspaces $V^1_i$,
it acts as  multiplication by $\la$ when $i\in (b_-,b_+)$,
and $0$ when $i\in [b_-+b_+,n]$.

Restricted to $V^2_i$, $i\in Y$, an RE  matrix $A$ of  
of Type 2 is equal to 
$
\begin{array}{||cc||}
\la& y_i\\
0        &0
\end{array}
$
if $i\in Y_-$
and
$
\begin{array}{||cc||}
\la& 0\\
y_i  &0
\end{array}
$
if $i\in Y_+$.
In either cases, it has eigenvalues $\la$ and $0$.
Observe that $Y_-\cup\si(Y_+)\subset [1,b]$ 
and $Y_-\cap\si(Y_+)=\emptyset $; thus
$\card\bigl(Y_-\cup\si(Y_+)\bigr)=\card(Y)=K$.
Therefore there are $b-K$ one-dimensional subspaces $V^1_i$
where $A$ acts as multiplication by $\la$.
So the eigenvalue $\la$ has multiplicity $b$. The 
zero eigenvalue enters with multiplicity $n-b$.
\end{proof}
Note that when one of the 
eigenvalues $\la,\mu$ tends to zero, the RE matrices
of Type 1 turns into  RE matrices  of Type 2.
Non-degenerate numerical RE matrices belong
to the first class. This is the case when $b_- + b_+-1=n$
and $\la\not =0\not =\mu$.
The RE matrices from \cite{KSS}
written out in Introduction are of this kind, and
the Type 1 solutions from Theorem \ref{clth} are their 
generalization to higher dimensions. Let us illustrate 
Theorem \ref{clth}
on those examples. 
\begin{enumerate}
\item
The matrix $D_n$ is a solution of Type 1.
Here, $\la=-\mu=y_i$,  $i\in Y$, where the set $Y$ 
coincides with $I$ if $n=0\!\! \mod 2$ and 
$Y=I\backslash\{\frac{n+1}{2}\}$ if $n=1\!\! \mod 2$.
The similarity transformation by a diagonal matrix
changes the elements $y_i$ however preserving the condition 
$y_i y_{n+1-i} =\la^2$.
\item
The matrix $A^{1,1}$ has $b_-=1$,
$b_+=2$, while for the matrix $A^{2,2}$ these 
parameters take the values  $b_-=2$ and $b_+=3$. These 
solutions have 
in common $(b_-, b_+)=\emptyset$. On the contrary,
the interval $(b_-, b_+)$ is not empty in the case of
matrices $A^{2,1}$ and  $A^{3,1}$. For them,
one has $b_-=1$,  $b_+=3$ and  $b_-=1$,  $b_+=4$.
\item
The diagonal matrix $P_k=\la\sum_{i=1}^k e^i_i$ is an RE matrix of 
Type 2 with $Y= \emptyset$  and $k=b\leq n$.
\end{enumerate}
\begin{remark}
Numerical RE matrices play an important role in
$\U_q\bigl(gl(n)\bigr)$-covariant 
quantization on the adjoint orbits in $\End^{\tp 2}(\C^n)$,
\cite{DoM1,DoM2}. They allow to explicitly represent
the quantized algebras of functions on orbits as quotients
of the RE algebra and, simultaneously, as subalgebras of functions on 
the quantum group. Theorem \ref{clth} classify all 
such realizations. It implies, in particular, that
all the symmetric and bisymmetric orbits  
(consisting of matrices with two and three different 
eigenvalues) can be represented in this way.
\end{remark}

\end{document}